\documentclass[11pt]{article}
\usepackage[active]{srcltx}
\usepackage{amssymb,amsmath}
\usepackage{graphicx}
\usepackage{amsbsy}
\usepackage{mathrsfs}

\newtheorem{proposition}{Proposition}[section]
\newtheorem{theorem}{Theorem}[section]
\newtheorem{lemma}{Lemma}[section]

\newtheorem{remark}{Remark}[section]

\numberwithin{equation}{section}

\begin{document}

\begin{center}{\large\sc
Rellich Inequalities on Finsler-Hadamard Manifolds}\\
\vspace{0.5cm} {\large Alexandru Krist\'aly \& Du\v san Repov\v s

}
\end{center}

\begin{abstract}
{\footnotesize \noindent In this paper we are dealing with improved
Rellich inequalities on Finsler-Hadamard manifolds with vanishing
mean covariation where the remainder terms are expressed by means of
the flag curvature. By exploiting various arguments from Finsler
geometry we show that more weighty curvature implies more powerful
improvements. The sharpness of the involved constants are also
studied. 
}
\end{abstract}



\vspace{0.0cm}

 \section{Introduction and main results}
The Hardy inequality
$$
    \int_{\mathbb R^n}|\nabla u|^2{\text d}x \geq \frac{(n-2)^2}{4}\int_{\mathbb
    R^n}\frac{u^2}{|x|^2}{\text d}x,\ \forall u\in C_0^\infty(\mathbb R^n),
$$
 plays a central role in the study of singular elliptic problems, $n\geq
 3$, where the constant $\frac{(n-2)^2}{4}$ is sharp but not
 achieved. The second-order Hardy inequalities are referred as {\it Rellich
 inequalities} whose most familiar forms can be stated as follows; given  $n\geq
 5$, one has
  \begin{equation}\label{standard-rellich-1}
\int_{\mathbb R^n}  (\Delta u)^2{\text d}x\geq
\frac{n^2(n-4)^2}{16}\int_{\mathbb R^n} \frac{u^2}{|x|^{4} }{\text
d}x, \forall u\in C_0^\infty(\mathbb R^n),
\end{equation}
\begin{equation}\label{standard-rellich-2}
\int_{\mathbb R^n} (\Delta u)^2{\text d}x\geq \frac{n^2}{4}
\int_{\mathbb R^n} \frac{|\nabla u|^2}{|x|^2}{\text d}x, \ \forall
u\in C_0^\infty(\mathbb R^n),
\end{equation}
where both constants $\frac{n^2(n-4)^2}{16}$ and $\frac{n^2}{4}$ are
sharp, but are never achieved. Hereafter, $\Delta$, $\nabla$,
$|\cdot|$ and ${\text d}x$ denote the classical Laplace operator,
the Euclidean gradient, the Euclidean norm and the Lebesgue measure
on $\mathbb R^n$, respectively. Due to the lack of extremal
functions in the Rellich inequalities, various improvements of
(\ref{standard-rellich-1}) and (\ref{standard-rellich-2}) can be
found in the literature; see e.g. Ghoussoub and  Moradifam
\cite{GM}, Tertikas and Zographopoulos \cite{TZ}, and references
therein.

 Hardy and Rellich inequalities have also been
studied on {\it curved spaces}. As far as we know,  Carron
\cite{Carron-JMPA} first studied Hardy inequalities on complete,
non-compact Riemannian manifolds. Motivated by \cite{Carron-JMPA},
Kombe and \"Ozaydin \cite{KO-TAMS-2009,KO-TAMS-2013}, and Yang, Su
and Kong \cite{YSK-Kinaiak} presented various Brezis-Vazquez-type
improvements of Hardy and Rellich inequalities on complete,
non-compact Riemannian manifolds. Recently, Krist\'aly
\cite{Kristaly} proved Hardy inequalities on reversible Finsler
manifolds where the improvements are given in terms of the
curvature.

The purpose of our paper is to describe improved Rellich
inequalities on {\it Finsler-Hadamard manifolds} (i.e., complete,
simply connected Finsler manifolds with non-positive flag curvature)
where the remainder terms involve the flag curvature. Two facts
should be highlighted:
\begin{itemize}
  \item  We prove that Rellich inequalities on
Finsler-Hadamard manifolds are better improved once the flag
curvature is more powerful. These phenomena can be considered as
second-order versions of the result described  in \cite{Kristaly}.
  \item Since Rellich inequalities on Finsler manifolds involve the
  highly nonlinear Finsler-Laplace operator $\boldsymbol{\Delta}$, expected properties  usually fail (which trivially hold on the 'linear' Riemannian
  context). Although our results are also genuinely new  in the Riemannian framework, we prefer to present them in
the context of Finsler geometry. In this manner, we emphasize the
deep connection between geometric and analytic phenomena which are
behind of second-order Sobolev-type inequalities on Finsler
manifolds, providing a new bridge between Finsler geometry and PDEs.
This fact is interesting in its own right as well from the point of
view of applications, see Antonelli, Ingarden and Matsumoto
\cite{AIM}.
\end{itemize}

In order to present the nature of our results, we  need some
notations and notions, see \S \ref{sect-2}.

Let $(M,F)$ be an $n-$dimensional complete reversible Finsler
manifold $(n\geq 5)$, $d_F:M\times M\to \mathbb R$ being the natural
distance function generated by the Finsler metric $F$, and let
$F^*:T^*M\to [0,\infty)$ be the polar transform of $F$. Let $D
u(x)\in T_x^*M$, $\boldsymbol{\nabla}u(x)\in T_xM$ and
$\boldsymbol{\Delta}u(x)$ be the derivative, gradient and
Finsler-Laplace operator of $u$ at $x\in M$, respectively.
 Let  ${\text d}V_F(x)$ be
 the Busemann-Hausdorff measure on $(M,F)$ and for a fixed $x_0\in M$, let us denote
$d(x):=d_F(x_0,x).$

Let $G_F:C_0^\infty(M)\to \mathbb R$ be defined by
$$G_F(u)=\int_M  \left[u(x)^2\boldsymbol{\Delta}(d(x)^{-2})-d(x)^{-2}\boldsymbol{\Delta}(u(x)^2)\right]{\text
d}V_F(x),$$ which gives the 'Green-deflection' of $u$ with respect
to the Finsler metric $F$; for a generic Finsler manifold $(M,F)$,
the function  $G_F$ does not vanish. However, $G_F\equiv 0$ whenever
$(M,F)$ is Riemannian due to Green's identity. Finally, we introduce
the following class of functions
$$C_{0,F}^\infty(M)=\left\{u\in C_0^\infty(M):G_F(u)=0\right\}.$$

A simple consequence of our main results (see Theorems \ref{main-1}
\& \ref{main-2}) can be stated as follows.

\begin{theorem}\label{Rellich-intro-1}
Let $(M,F)$ be an $n-$dimensional reversible Finsler-Hadamard
manifold with vanishing mean covariation, and suppose the flag
curvature on $(M,F)$ is bounded above by ${c}\leq 0$.
\begin{itemize}
  \item[{\rm (a)}] If $n\geq 5$, then for every $u\in C_{0,F}^\infty(M)$ one has  \begin{eqnarray*}
 \int_M (\boldsymbol{\Delta} u)^2\text{\emph{d}}V_F(x)  &\geq& \frac{n^2(n-4)^2}{16}\int_M
\frac{u^2}{d(x)^{4}}\text{\emph{d}}V_F(x) \\
   &&+\frac{3|{c}|n(n-1)(n-2)(n-4)}{4}\int_M
\frac{u^2}{(\pi^2+|{c}|d(x)^{2})d(x)^{2}}\text{\emph{d}}V_F(x),
\end{eqnarray*}
and the constant $\frac{n^2(n-4)^2}{16}$ is sharp.
  \item[{\rm (b)}] If $n\geq 9$, then for every $u\in  C_{0,F}^\infty(M)$ one has
  \begin{eqnarray*}
 \int_M (\boldsymbol{\Delta} u)^2\text{\emph{d}}V_F(x)  &\geq& \frac{n^2}{4}\int_M
\frac{F^*(x,Du(x))^2}{d(x)^{2}}\text{\emph{d}}V_F(x), \\&&
   +\frac{3|{c}|n(n-1)(n-4)^2}{8}\int_M
\frac{u^2}{(\pi^2+|{c}|d(x)^{2})d(x)^{2}}\text{\emph{d}}V_F(x),
\end{eqnarray*}
and the constant $\frac{n^2}{4}$ is sharp.
\end{itemize}
\end{theorem}

\begin{remark}\rm
(i) When the flag curvature on $(M,F)$ becomes more powerful (i.e.,
$|{c}|$ is large), the Rellich inequalities in Theorem
\ref{Rellich-intro-1} is also better improved.

(ii) Theorem \ref{Rellich-intro-1} is also new for Cartan-type
Riemannian manifolds; indeed, these spaces belong to the class of
Cartan-Finsler manifolds with vanishing mean covariation and
$C_{0,F}^\infty(M)=C_0^\infty(M)$.
\end{remark}

In Section \ref{sect-2} we shall recall some elements from Finsler
geometry, namely  the flag curvature, Laplace and volume
comparisons, differentials. In Section \ref{sect-3} we shall prove
our main results (see Theorems \ref{main-1} \& \ref{main-2}), while
in Section \ref{sect-concluding} we shall present some concluding
remarks.

\section{Preliminaries}\label{sect-2}


 Let $M$ be a connected
$n-$dimensional $C^{\infty}$ manifold and $TM=\bigcup_{x \in M}T_{x}
M $ its tangent bundle. The pair $(M,F)$ is called a reversible
\textit{Finsler manifold} if the continuous function $F:TM\to
[0,\infty)$ satisfies the following conditions

(a) $F\in C^{\infty}(TM\setminus\{ 0 \});$

(b) $F(x,ty)=|t|F(x,y)$ for all $t\in \mathbb R$ and $(x,y)\in TM;$

(c) $g_{ij}(x,y):=[\frac12F^{2}%
(x,y)]_{y^{i}y^{j}}$ is positive definite for all $(x,y)\in
TM\setminus\{ 0 \}.$

 \noindent
 If $g_{ij}(x)=g_{ij}(x,y)$ is independent of $y$
then $(M,F)$ is called {\it Riemannian manifold}.  A {\it Minkowski
space} consists of a finite dimensional vector space $V$ and a
Minkowski norm which induces a Finsler metric on $V$ by translation,
i.e., $F(x,y)$ is independent of the base point $x$; in such cases
we often write $F(y)$ instead of $F(x,y)$. While there is a unique
Euclidean space (up to isometry), there are infinitely many
(isometrically different) Minkowski spaces.

We consider the {\it polar transform} of $F$, defined for every
$(x,\xi)\in T^*M$ by
\begin{equation}\label{polar-transform}
    F^*(x,\xi)=\sup_{y\in T_xM\setminus
    \{0\}}\frac{\xi(y)}{F(x,y)}.
\end{equation}
Note that for every $x\in M$, the function $F^*(x,\cdot)$ is a
Minkowski norm on $T_x^*M.$ Since $F^*(x,\cdot)^2$ is twice
differentiable on $T_x^*M\setminus \{0\}$, we consider the matrix $$g_{ij}^*(x,\xi):=[\frac12F^*(x,\xi)^2%
]_{\xi^{i}\xi^{j}}$$ for every $\xi=\sum_{i=1}^n\xi^i{\rm d}x^i\in
T_x^*M\setminus \{0\}$  in a local coordinate system $(x^i).$

Let $\pi ^{*}TM$ be the pull-back bundle of the tangent bundle $TM$
generated by the natural projection $\pi:TM\setminus\{ 0 \}\to M,$
see Bao, Chern and Shen \cite[p. 28]{BCS}. The vectors of the
pull-back bundle $\pi ^{*}TM$ are denoted by $(v;w)$ with
$(x,y)=v\in TM\setminus\{ 0 \}$ and $w\in T_xM.$ For simplicity, let
$\partial_i|_v=(v;\partial/\partial x^i|_x)$ be the natural local
basis for $\pi ^{*}TM$, where $v\in T_xM.$ One can introduce the
{\it fundamental tensor} $g$ on $\pi ^{*}TM$  by
\begin{equation}\label{funda-Cartan-tensors}
    g^v:=g(\partial_i|_v,\partial_i|_v)=g_{ij}(x,y),
\end{equation}
where $v=y^i{(\partial}/{\partial x^i})|_x.$  Unlike the Levi-Civita
connection \index{Levi-Civita connection} in the Riemannian case,
there is no unique natural connection in the Finsler geometry. Among
all natural connections on the pull-back bundle $\pi ^{*}TM,$ we
choose a torsion free and almost metric-compatible linear connection
on $\pi ^{*}TM$, the so-called \textit{Chern connection}, see Bao,
Chern and Shen \cite[Theorem 2.4.1]{BCS}.  The coefficients of the
Chern connection are denoted by $\Gamma_{jk}^{i}$, which replace the
well known Christoffel symbols from Riemannian geometry. A Finsler
manifold is said to be of \textit{Berwald type} if the coefficients
$\Gamma_{ij}^{k}(x,y)$ in natural coordinates are independent of
$y$. It is clear that Riemannian manifolds and $($locally$)$
Minkowski spaces are Berwald spaces. The Chern connection induces in
a natural manner on $\pi ^{*}TM$ the \textit{curvature tensor} $R$,
see Bao, Chern and Shen \cite[Chapter 3]{BCS}. By means of the
connection, we also have
 the {\it covariant derivative} $D_vu$ of a vector field
$u$ in the direction $v\in T_xM.$  Note that $v\mapsto D_vu$ is not
linear. A vector field $u=u(t)$ along a curve $\sigma$ is said to be
{\it parallel} if $D_{\dot \sigma}u=0.$ A $C^\infty$ curve
$\sigma:[0,a]\to M$ is called a {\it geodesic} if  $D_{\dot
\sigma}{\dot \sigma}=0.$ Geodesics are considered to be parametrized
proportionally to their arc-length. The Finsler manifold is said to
be {\it complete} if every geodesic segment can be extended to
$\mathbb R.$

 Let $u,v\in T_xM$ be two non-collinear vectors and
$\mathcal S={\rm span}\{u,v\}\subset T_xM$. By means of the
curvature tensor $R$, the {\it flag curvature} of the flag
$\{\mathcal S,v\}$ is then defined by
\begin{equation} \label{ref-flag}
K(\mathcal S;v) =\frac{g^v(R(U,V)V, U)}{g^v(V,V) g^v(U,U) -
g^v(U,V)^{2}},
\end{equation}
where $U=(v;u),V=(v;v)\in \pi^*TM.$ If for some ${c}\in \mathbb R$
we have  $K(\mathcal S;v)\leq {c}$  for every choice of $U$ and $V$,
we say that the flag curvature is bounded from above by ${c}$ and we
write ${\bf K}\leq {c}$.  $(M,F)$ is called a {\it Finsler-Hadamard}
manifold if it is complete, simply connected and ${\bf K}\leq 0.$ If
$(M,F)$ is Riemannian, the flag curvature reduces to the well known
sectional curvature.

  Let $\sigma: [0,r]\to
M$ be a piecewise $C^{\infty}$ curve. The value $ L_F(\sigma)=
\int_{0}^{r} F(\sigma(t), \dot\sigma(t))\,{\text d}t $ denotes the
\textit{integral length} of $\sigma.$  For $x_1,x_2\in M$, denote by
$\Lambda(x_1,x_2)$ the set of all piecewise $C^{\infty}$ curves
$\sigma:[0,r]\to M$ such that $\sigma(0)=x_1$ and $\sigma(r)=x_2$.
Define the {\it distance function} $d_{F}: M\times M \to[0,\infty)$
by
\begin{equation}\label{quasi-metric}
  d_{F}(x_1,x_2) = \inf_{\sigma\in\Lambda(x_1,x_2)}
 L_F(\sigma).
\end{equation}
Clearly, $d_F$ satisfies all properties of the metric (i.e.,
$d_{F}(x_1,x_2) =0$ if and only if $x_1=x_2,$ $d_F$ is symmetric,
and it satisfies the triangle inequality). The open {\it metric
ball} with center $x_0\in M$ and radius $\rho>0$ is defined by
$B(x_0,\rho)=\{x\in M:d_F(x_0,x)< \rho\}.$

 Let $\{{\partial}/{\partial x^i}
\}_{i=1,...,n}$ be a local basis for the tangent bundle $TM,$ and
let $\{{\rm d}x^i \}_{i=1,...,n}$ be its dual basis for $T^*M.$ Let
$B_x(1)=\{y=(y^i):F(x,y^i \partial/\partial x^i)< 1\}$ be the unit
tangent ball at $T_xM$.  The {\it Busemann-Hausdorff  volume form}
${\text d}V_F$ on $(M,F)$ is defined by
\begin{equation}\label{volume-form}
    {\text d}V_F(x)=\sigma_F(x){\text d}x^1\wedge...\wedge {\text d}x^n,
\end{equation}
where $\sigma_F(x)=\frac{\omega_n}{{\rm Vol}(B_x(1))}$. Hereafter,
$\omega_n$ will denote the volume of the unit $n-$dimensional ball
and Vol$(S)$ the Euclidean volume of the set $S\subset \mathbb R^n$.
The {\it Finslerian-volume} of a bounded open set $S\subset M$ is
defined as Vol$_F(S)=\int_S {\text d}V_F(x)$. In general, one has
that for every $x\in M,$
\begin{equation}\label{volume-comp-nullaban}
   \lim_{\rho\to 0^+}\frac{{\rm Vol}_F(B(x,\rho))}{\omega_n
   \rho^n}=1.
\end{equation} When $(\mathbb
R^n,F)$ is a Minkowski space, then by virtue of (\ref{volume-form}),
Vol$_F(B(x,\rho))=\omega_n\rho ^n$ for every $\rho>0$ and $x\in
\mathbb R^n.$

  The {\it Legendre transform}
$J^*:T^*M\to TM$ associates to each element $\xi\in T_x^*M$ the
unique maximizer on $T_xM$ of the map $y\mapsto
\xi(y)-\frac{1}{2}F^2(x,y)$. This element can also be interpreted as
the unique vector $y\in T_xM$ with the following properties
\begin{equation}\label{ohta-sturm-dual}
    F(x,y)=F^*(x,\xi)\ {\rm and}\ \xi(y)=F(x,y)F^*(x,\xi).
\end{equation}
In a similar manner we can define the Legendre transform $J:TM\to
T^*M$. In particular, $J^*=J^{-1}$ on $T_x^*M$ and if
$\xi=\sum_{i=1}^n\xi^i {\rm d}x^i\in T_x^*M$ and $y=\sum_{i=1}^ny^i
(\partial/\partial x^i)\in T_xM$, then one has
\begin{equation}\label{j-csillag}
    J(x,y)=\sum_{i=1}^n\frac{\partial}{\partial
y^i}\left(\frac{1}{2}F(x,y)^2\right)\frac{\partial}{\partial x^i}\ \
{\rm and} \ \ J^*(x,\xi)=\sum_{i=1}^n\frac{\partial}{\partial
\xi^i}\left(\frac{1}{2}F^*(x,\xi)^2\right)\frac{\partial}{\partial
x^i}.
\end{equation}
  Let $u:M\to \mathbb R$ be a differentiable function in the
distributional sense. The {\it gradient} of $u$ is defined by
\begin{equation}\label{gradient-fins}
    \boldsymbol{\nabla} u(x)=J^*(x,Du(x)),
\end{equation}
where $Du(x)\in T_x^*M$ denotes the (distributional) {\it
derivative} of $u$ at $x\in M.$ In general,
 $u\mapsto\boldsymbol{\nabla} u $ is not linear.

 Let $x_0\in M$ be fixed. From now on when no confusion arises, we shall introduce the abbreviation
\begin{equation}\label{abbreviation}
    d(x)=d_F(x_0,x).
\end{equation}
 Due to Ohta and Sturm
 \cite{Ohta-Sturm} and by relation (\ref{ohta-sturm-dual}), one has
{ \begin{equation}\label{tavolsag-derivalt}
    F(x,\boldsymbol{\nabla} d(x))=F^*(x,D d(x))=D d(x)(\boldsymbol{\nabla} d(x))=1\ {\rm for\ a.e.}\ x\in
    M.
\end{equation}}
In fact, relations from (\ref{tavolsag-derivalt}) are valid for
every  $x\in M\setminus (\{x_0\}\cup{\rm Cut}(x_0))$, where ${\rm
Cut}(x_0)$ denotes the cut locus of $x_0,$ see Bao, Chern and Shen
\cite[Chapter 8]{BCS}. Note that ${\rm Cut}(x_0)$ has null Lebesgue
(thus Hausdorff) measure for every $x_0\in M$.

Let $X$ be a vector field on $M$. In a local coordinate system
$(x^i)$, by virtue of (\ref{volume-form}), the {\it divergence} is
defined by div$(X)=\frac{1}{\sigma_F}\frac{\partial}{\partial
x^i}(\sigma_F X^i).$ The {\it Finsler-Laplace operator}
 $$\boldsymbol{\Delta} u={\rm div}(\boldsymbol{\nabla} u)$$ acts on $W^{1,2}_{\rm
 loc}(M)$ and for every $v\in C_0^\infty(M)$, we have
\begin{equation}\label{Green}
\int_M v\boldsymbol{\Delta} u {\text d}V_F(x)=-\int_M
Dv(\boldsymbol{\nabla} u){\text d}V_F(x),
\end{equation}
see  Ohta and Sturm
 \cite{Ohta-Sturm} and Shen \cite{Shen-monograph}. In the Riemannian
 case, the Finsler-Laplace operator reduces to the Laplace-Beltrami operator, see Bonanno, G. Molica Bisci, V. R\u
 adulescu \cite{BMBR}.

 Let $\{e_i\}_{i=1,...,n}$ be a
basis for $T_xM$ and $g_{ij}^v=g^v(e_i,e_j)$. The {\it mean
distortion} $\mu:TM\setminus \{0\}\to (0,\infty)$ is defined by
$\mu(v)=\frac{\sqrt{{\rm det}(g_{ij}^v)}}{\sigma_F}$. The {\it mean
covariation} ${\bf S}:TM\setminus\{0\}\to \mathbb R$ is defined by
$${\bf S}(x,v)=\frac{d}{dt}(\ln \mu(\dot\sigma_v(t)))\big|_{t=0},$$ where
$\sigma_v$ is the geodesic such that $\sigma_v(0)=x$ and $\dot
\sigma_v(0)=v.$ 
We say that $(M,F)$ has {\it vanishing mean covariation} if ${\bf
S}(x,v)= 0$ for every $(x,v)\in TM$, and we denote this by ${\bf
S}=0$. We recall that any Berwald space has vanishing mean
covariation, see Shen \cite{Shen-volume}.

 We conclude this section by some important comparison results.  Let $x_0\in M$ be fixed and
recall the notation introduced in (\ref{abbreviation}). First, one
has
\begin{equation}\label{local-laplace}
\boldsymbol{\Delta}d(x)-\frac{n-1}{d(x)}=o(1)\ {\rm as}\ x\to x_0.
\end{equation}
In order to have a global estimate for $\boldsymbol{\Delta}d(x)$, we
consider for every ${c}\leq 0$ the function ${\bf
ct}_{c}:(0,\infty)\to \mathbb R$ defined by
$${\bf ct}_{c}(\rho)=\left\{
  \begin{array}{lll}
    \frac{1}{\rho}
    & \hbox{if} &  {c}=0, \\
    \\
  \sqrt{|{c}|}\coth(\sqrt{|{c}|}\rho) & \hbox{if} & {c}<0.
  \end{array}\right.$$

\begin{theorem}\label{comparison-laplace} Let $(M,F)$ be an $n-$di\-men\-sional Finsler-Hadamard manifold with
${\bf S}=0$ and ${\bf K}\leq {c}\leq 0$, and let $x_0\in M$ be
fixed. Then the following assertions hold:

\begin{itemize}
\item[{\rm (a)}] {\rm (see \cite[Theorem 5.1]{Wu-Xin})} For a.e. $x\in
    M$ one has    $\boldsymbol{\Delta}d(x)\geq (n-1){\bf ct}_{c}(d(x)).$
  \item[{\rm (b)}] {\rm (see \cite[Theorem 6.1]{Wu-Xin})} The function
$\rho\mapsto \frac{{\rm Vol}_F(B(x,\rho))}{\rho^n}$ is
non-decreasing, $\rho>0$. In particular, by {\rm
(\ref{volume-comp-nullaban})} we have
$$
{{\rm Vol}_F(B(x,\rho))}\geq \omega_n \rho^n\ {for\ all}\ x\in M\
{and}\ \rho>0.
$$
\end{itemize}
\end{theorem}

\section{Main results}\label{sect-3}

Let ${\bf D}_c:[0,\infty)\to \mathbb R$ be the function defined by
$${\bf D}_{c}(\rho)=\left\{
  \begin{array}{lll}
    0
    & \hbox{if} &  \rho=0, \\
    \\
  \rho {\bf ct}_c(\rho)-1 & \hbox{if} & \rho>0.
  \end{array}\right.$$
 It is clear that ${\bf D}_c\geq 0.$

 In order to establish our main results, we first need a quantitative
 Hardy inequality; see \cite{Kristaly} for a particular form.  For the reader's convenience we provide its proof.

\begin{lemma}\label{lemma-Hardy}  Let $(M,F)$ be an $n-$di\-men\-sional Finsler-Hadamard manifold with
${\bf S}=0$ and let ${\bf K}\leq {c}\leq 0$,  $x_0\in M$ be fixed,
and choose any $\alpha\in \mathbb R$ such that $n-2+\alpha>0.$
 Then for every $u\in C_0^\infty(M)$ we have
{ \begin{eqnarray*}
  \int_{M}d(x)^\alpha F^*(x,Du(x))^2\text{\emph{d}}V_F(x) &\geq &
  \frac{(n-2+\alpha)^2}{4}\int_{M}{d(x)^{\alpha-2}}{u(x)^2}\text{\emph{d}}V_F(x)\\&&+\frac{(n-2+\alpha)(n-1)}{2}\int_{M}{d(x)^{\alpha-2}}{\bf
D}_c(d(x)){u(x)^2}\text{\emph{d}}V_F(x).
\end{eqnarray*}}
\end{lemma}

{\it Proof.} By convexity  and (\ref{j-csillag}), one has
\begin{equation}\label{f-cillag-Taylor}
    F^*(x,\xi_2)^2\geq F^*(x,\xi_1)^2+2(\xi_2-\xi_1)(J^*(x,\xi_1)), \
\forall \xi_1,\xi_2\in T_x^*M.
\end{equation}
Let $u\in C_0^\infty(M)$ be arbitrarily and choose
$\tau=\frac{n-2+\alpha}{2}>0$. Let $v(x)=d(x)^\tau u(x)$. Therefore,
for $u(x)=d(x)^{-\tau} v(x)$ one has $Du(x)=-\tau d(x)^{-\tau-1}v(x)
Dd(x)+d(x)^{-\tau}Dv(x).$ By inequality (\ref{f-cillag-Taylor})
applied for $\xi_2=-Du(x)$ and $\xi_1=\tau d(x)^{-\tau-1}v(x)
Dd(x)$, the symmetry of $F^*(x,\cdot)$ implies that
\begin{eqnarray*}
 F^*(x,Du(x))^2  &=& F^*(x,-Du(x))^2 \\
   &\geq& F^*(x,\tau d(x)^{-\tau-1}v(x) Dd(x))^2 -2d(x)^{-\tau}Dv(x)(J^*(x,\tau d(x)^{-\tau-1}v(x) Dd(x))).
\end{eqnarray*}
Since $F^*(x,Dd(x))=1$ (see (\ref{tavolsag-derivalt})),
$J^*(x,Dd(x))=\boldsymbol{\nabla}d(x)$ and $Dv(x)\in T_x^*M,$ we
obtain
$$F^*(x,Du(x))^2\geq \tau^2 d(x)^{-2\tau-2}v(x)^2-2\tau  d(x)^{-2\tau-1}v(x)Dv(x)(\boldsymbol{\nabla}d(x)).$$
Multiplying the latter inequality by $d(x)^\alpha$, and integrating
 over $M$, we obtain
$$\int_{M}d(x)^\alpha F^*(x,Du(x))^2{\text d}V_F(x)  \geq  \tau^2\int_{M}d(x)^{\alpha-2\tau-2}v(x)^2{\text d}V_F(x)+R_0,$$
where
\begin{eqnarray*}
  R_0 &=& -2\tau\int_M
d(x)^{\alpha-2\tau-1}v(x)Dv(x)(\boldsymbol{\nabla}d(x)){\text d}V_F(x) \\
  &=&-\frac{\tau}{\alpha-2\tau} \int_M D(v(x)^2)(\boldsymbol{\nabla}(d(x)^{\alpha-2\tau})){\text d}V_F(x)\\
  &=&\frac{\tau}{\alpha-2\tau} \int_M
v(x)^2\boldsymbol{\Delta}(d(x)^{\alpha-2\tau}){\text d}V_F(x) \ \ \
\ \ \ \ \ \ \ \ \ \ \ \ \ \ \ \ \ \ \ \ \ \ \ \ \ \ \ \ \ \ \ \ \ \
\ \ \ \ \ \ \ {\rm
(see\ (\ref{Green}))}\\
&=& \tau \int_M
u(x)^2d(x)^{\alpha-2}\left[\alpha-2\tau-1+d(x)\boldsymbol{\Delta}d(x)\right]{\text d}V_F(x)\\
&\geq & \tau(n-1)\int_M u(x)^2d(x)^{\alpha-2}\left[d(x){\bf
ct}_c(d(x))-1\right] {\text d}V_F(x), \ \ \
\ \ \ \ \ \ \ \ {\rm (see\ Theorem\ \ref{comparison-laplace}\ (a))}\\
&= & \tau(n-1)\int_M d(x)^{\alpha-2}{\bf D}_c(d(x))u(x)^2 {\text
d}V_F(x),
\end{eqnarray*}
which completes the  proof. \hfill $\square$\\

For every $x\in M$ and $y\in T_xM$, $\xi\in T_x^*M,$ we introduce
the function
\begin{equation}\label{K-function}
    K_F(x,y,\xi)=\xi(y)-J(x,y)(J^*(x,\xi)).
\end{equation}
For $\alpha\in \mathbb R$ with $n-4+\alpha>0$ we introduce the {\it
Green-deflection} function
 $G_F^\alpha:C_0^\infty(M)\to \mathbb R$ defined by
$$G_F^\alpha(u)=\int_M  K_F\left(x,\boldsymbol{\nabla}(u(x))^2,D(d(x)^{\alpha-2})\right){\text
d}V_F(x).$$ The layer cake representation and the fact that
$n-4+\alpha>0$ imply that the function $G_F^\alpha$ is well defined.
Moreover, by  definition of $K_F$ and relations
(\ref{gradient-fins}) and (\ref{Green}) one has
\begin{equation}\label{g_F-repres}
    G_F^\alpha(u)=\int_M
\left[u(x)^2\boldsymbol{\Delta}(d(x)^{\alpha-2})-d(x)^{\alpha-2}\boldsymbol{\Delta}(u(x)^2)\right]{\text
d}V_F(x).
\end{equation}
It is now clear that $G_F^\alpha\equiv 0$ whenever $(M,F)$ is
Riemannian due to Green's identity. In fact, the latter statement
also holds by the following observation.

\begin{proposition}\label{prop-eukl-green} $K_F\equiv 0$ if and only
if $(M,F)$ is Riemannian.
\end{proposition}

 {\it Proof.}  If $(M,F)$ is Riemannian then $g(x,y)=a(x)$, where
 $a(x)$ is a symmetric and positive-definite matrix and by Riesz representation, one can identify $T_xM$ and $T_x^*M.$ Moreover,
 $J(x,y)=a(x)y$ and $J^*(x,\xi)=a(x)^{-1}\xi.$ Consequently, we have
 $$K_F(x,y,\xi)=\xi(y)-J(x,y)(J^*(x,\xi))=\xi(y)-a(x)y(a(x)^{-1}\xi)=0.$$

Conversely, we assume that  $K_F\equiv 0$, i.e.,
$\xi(y)-J(x,y)(J^*(x,\xi))=0$ for every $x\in M$, $y\in T_xM$ and
$\xi\in T_x^*M.$ For an arbitrary $z\in T_xM$ replace $\xi=J(x,z)\in
T_x^*M$ into the preceding relation to obtain $J(x,z)(y)=J(x,y)(z)$.
In particular, $J(x,\cdot)$ is linear; by virtue of
(\ref{j-csillag}) it implies that $F(x,\cdot)^2$ comes from an inner
product on $T_xM.$
 \hfill $\square$ \\

Let us consider the following set of functions
$$C_{0,F,\alpha}^\infty(M)=\left\{u\in
C_0^\infty(M):G_F^\alpha(u)=0\right\}.$$ By Proposition
\ref{prop-eukl-green}, $C_{0,F,\alpha}^\infty(M)=C_{0}^\infty(M)$
whenever $(M,F)$ is Riemannian. However, in the generic Finsler
context the role of $C_{0,F,\alpha}^\infty(M)$ seems to be
indispensable for the study of Rellich inequalities.

We are in  position to state our first main result.

\begin{theorem}\label{main-1} {\rm (Rellich inequality I)} Let $(M,F)$ be an $n-$di\-men\-sional Finsler-Hadamard manifold with
${\bf S}=0$ and ${\bf K}\leq {c}\leq 0$, let $x_0\in M$ be fixed,
and choose any $\alpha\in \mathbb R$ such that $n-4+\alpha>0$ and
$\alpha<2.$
 Then for every $u\in C_{0,F,\alpha}^\infty(M)$ we have
{ \begin{eqnarray*}
  \int_{M}d(x)^\alpha (\boldsymbol{\Delta}u(x))^2\text{\emph{d}}V_F(x) &\geq &
  \frac{(n-4+\alpha)^2(n-\alpha)^2}{16}\int_{M}{d(x)^{\alpha-4}}{u(x)^2}\text{\emph{d}}V_F(x)\\&&+\frac{(n-4+\alpha)(n-\alpha)(n-2)(n-1)}{4}\int_{M}{d(x)^{\alpha-4}}{\bf
D}_c(d(x)){u(x)^2}\text{\emph{d}}V_F(x).
\end{eqnarray*}}
Moreover, the constant $\frac{(n-4+\alpha)^2(n-\alpha)^2}{16}$ is
sharp.
\end{theorem}
 {\it Proof.} Throughout  the proof, we shall consider $\gamma=\frac{n-4+\alpha}{2}>0.$ Since
$\alpha<2$, a simple calculation and Theorem
\ref{comparison-laplace}(a) yield
\begin{eqnarray*}
  \boldsymbol{\Delta}(d(x)^{\alpha-2}) &=& (\alpha-2)[\alpha-3+d(x)\boldsymbol{\Delta}(d(x))]d(x)^{\alpha-4} \\
   &\leq&
   (\alpha-2)[\alpha-3+(n-1)d(x){\bf ct}_c(d(x))]d(x)^{\alpha-4}\\
   &=&
   (\alpha-2)\left[2\gamma+(n-1){\bf D}_c(d(x))\right]d(x)^{\alpha-4}.
\end{eqnarray*}
Let us fix $u\in C_{0,F,\alpha}^\infty(M)$. Multiplying the above
inequality by $u^2$, we see that
\begin{equation}\label{re-1-1}
\int_{M}\boldsymbol{\Delta}(d(x)^{\alpha-2})u(x)^2{\text
d}V_F(x)\leq (\alpha-2)\int_M[2\gamma+(n-1){\bf
D}_c(d(x))]d(x)^{\alpha-4}u(x)^2{\text d}V_F(x).
\end{equation}
Note that
$$\boldsymbol{\Delta}(u(x)^2)=2{\rm div} (u\boldsymbol{\nabla}(u(x))=2F^*(x,Du(x))^2+2u\boldsymbol{\Delta}(u(x)).$$
Multiplying the latter relation by $d^{\alpha-2}$ and integrating
over $M$, we obtain
$$\int_Md(x)^{\alpha-2}\boldsymbol{\Delta}(u(x)^2){\text d}V_F(x)=2\int_Md(x)^{\alpha-2}F^*(x,Du(x))^2{\text d}V_F(x)+2\int_Md(x)^{\alpha-2}u\boldsymbol{\Delta}(u(x)){\text d}V_F(x).$$
Subtracting the latter relation by (\ref{re-1-1}), one gets that
\begin{eqnarray*}
  G_F^\alpha(u) &\leq& (\alpha-2)\int_M[2\gamma+(n-1){\bf
D}_c(d(x))]d(x)^{\alpha-4}u(x)^2{\text d}V_F(x) \\
   && -2\int_Md(x)^{\alpha-2}F^*(x,Du(x))^2{\text d}V_F(x)-2\int_Md(x)^{\alpha-2}u\boldsymbol{\Delta}(u(x)){\text d}V_F(x).
\end{eqnarray*}
Since $u\in C_{0,F,\alpha}^\infty(M)$, then $G_F^\alpha(u)=0$ and we
obtain that
\begin{eqnarray}\label{mindjart-kell}\nonumber -\int_Md(x)^{\alpha-2}u\boldsymbol{\Delta}(u(x)){\text d}V_F(x)&\geq&\frac{2-\alpha}{2}\int_M[2\gamma+(n-1){\bf
D}_c(d(x))]d(x)^{\alpha-4}u(x)^2{\text d}V_F(x)\\&&+
    \int_Md(x)^{\alpha-2}F^*(x,Du(x))^2{\text d}V_F(x).\end{eqnarray}
For the latter term we apply the Hardy inequality (Lemma
\ref{lemma-Hardy}),  and obtain {
\begin{eqnarray}\label{hardy-majd-a-vegen-kell}
  \int_{M}d(x)^{\alpha-2} F^*(x,Du(x))^2{\text d}V_F(x) &\geq &
 \gamma^2\int_{M}{d(x)^{\alpha-4}}{u(x)^2}{\text d}V_F(x)\nonumber\\&&+\gamma(n-1)\int_{M}{d(x)^{\alpha-4}}{\bf
D}_c(d(x)){u(x)^2}{\text d}V_F(x).
\end{eqnarray}}
Combining these inequalities, a trivial rearrangement now yields
\begin{eqnarray*} -\int_Md(x)^{\alpha-2}u\boldsymbol{\Delta}(u(x)){\text d}V_F(x)&\geq&\frac{\gamma(n-\alpha)}{2}\int_{M}{d(x)^{\alpha-4}}{u(x)^2}{\text d}V_F(x)\\&&+
    \frac{(n-1)(n-2)}{2}\int_{M}{d(x)^{\alpha-4}}{\bf
D}_c(d(x)){u(x)^2}{\text d}V_F(x).\end{eqnarray*}

The H\"older inequality for the left hand side of the above
inequality gives that
\begin{equation}\label{arith-mean}
     \left(\int_{M}d(x)^\alpha (\boldsymbol{\Delta}u(x))^2{\text d}V_F(x)\right)^\frac{1}{2}\cdot
\left(\int_{M}{d(x)^{\alpha-4}}{u(x)^2}{\text
d}V_F(x)\right)^\frac{1}{2}\geq
\int_Md(x)^{\alpha-2}|u\boldsymbol{\Delta}(u(x))|{\text d}V_F(x).
\end{equation}
The last inequalities and a simple estimate show that
\begin{eqnarray*}
  \int_{M}d(x)^\alpha (\boldsymbol{\Delta}u(x))^2{\text d}V_F(x) &\geq& \frac{\gamma^2(n-\alpha)^2}{4}\int_{M}{d(x)^{\alpha-4}}{u(x)^2}{\text d}V_F(x) \\
   &&+\frac{\gamma(n-\alpha)(n-2)(n-1)}{2}\int_{M}{d(x)^{\alpha-4}}{\bf
D}_c(d(x)){u(x)^2}{\text d}V_F(x),
\end{eqnarray*}
which completes the proof of Rellich inequality I.

 Now, we shall prove that in the Rellich inequality I  the constant
 $\tilde C:=\frac{\gamma^2(n-\alpha)^2}{4}$ is sharp. Clearly, it is enough to prove that
\begin{equation}\label{best-const}
    \tilde C=\inf_{u\in
    C_{0,F,\alpha}^\infty(M)\setminus \{0\}}\frac{\int_{M}d(x)^\alpha (\boldsymbol{\Delta}u(x))^2{\text d}V_F(x)}{\int_{M}{d(x)^{\alpha-4}}{u(x)^2}{\text
    d}V_F(x)}.
\end{equation}

 First, it follows by
 (\ref{local-laplace}) that there exists $0<r_0<\frac{n-\alpha}{2}$ such
 that $$\left|\boldsymbol{\Delta}d(x)-\frac{n-1}{d(x)}\right|\leq 1\ {\rm for\ a.e.}\ x\in B(x_0,r_0).$$
 In particular, one has
\begin{equation}\label{bx-0=ro}
    |-\gamma-1+d(x)\boldsymbol{\Delta}d(x)|\leq
    \frac{n-\alpha}{2}+d(x)\ \ \
{\rm for\ a.e.}\ x\in B(x_0,r_0).
\end{equation}

Let us fix numbers $r,R\in \mathbb R$ such that $0<r<R<r_0$ and a
smooth cutoff function $\psi:M\to [0,1]$ with
supp$(\psi)=\overline{B(x_0,R)}$ and $\psi(x)=1$ for $x\in
B(x_0,r).$ For every $0<\varepsilon<r,$ let
\begin{equation}\label{u-eps-mindenki-ezt}
    u_\varepsilon(x)=(\max\{\varepsilon,d(x)\})^{-\gamma},\ x\in
M.
\end{equation}
Note that $\psi u_\varepsilon$ can be approximated by elements from
$C_0^\infty(M)$ and since both functions $\psi$ and $u_\varepsilon$
are $d(x)-$radial, it follows by the representation
(\ref{g_F-repres}) of $G_F^\alpha$ that $G_F^\alpha(\psi
u_\varepsilon)=0,$ therefore, $\psi u_\varepsilon\in
C_{0,F,\alpha}^\infty(M)$ for every $0<\varepsilon<r.$

 One the one hand, by relation (\ref{bx-0=ro}) one has
\begin{eqnarray*}
  I_1(\varepsilon)&:=&\int_{M}d(x)^\alpha (\boldsymbol{\Delta}(\psi(x) u_\varepsilon(x)))^2{\text d}V_F(x)\\ &=&
  \int_{B(x_0,r)\setminus  B(x_0,\varepsilon)}d(x)^\alpha (\boldsymbol{\Delta}( d(x)^{-\gamma}))^2{\text d}V_F(x) \\&&+
  \int_{B(x_0,R)\setminus B(x_0,r)}d(x)^\alpha (\boldsymbol{\Delta}( \psi(x)d(x)^{-\gamma}))^2{\text d}V_F(x)\\
   &=&  \gamma^2\int_{B(x_0,r)\setminus  B(x_0,\varepsilon)}d(x)^{\alpha-2\gamma-4} [-\gamma-1+d(x)\boldsymbol{\Delta} d(x)]^2{\text
   d}V_F(x)+c(\alpha,r,R)\\&\leq&
   \gamma^2\left(\frac{n-\alpha}{2}+r\right)^2\int_{B(x_0,r)\setminus
   B(x_0,\varepsilon)}d(x)^{\alpha-2\gamma-4}{\text
   d}V_F(x)+c(\alpha,r,R)\\&=& \gamma^2\left(\frac{n-\alpha}{2}+r\right)^2\tilde
   I(\varepsilon)+c(\alpha,r,R),
\end{eqnarray*}
where
$$\tilde I(\varepsilon)=\int_{B(x_0,r)\setminus
   B(x_0,\varepsilon)}d(x)^{\alpha-2\gamma-4}{\text
   d}V_F(x)=\int_{B(x_0,r)\setminus
   B(x_0,\varepsilon)}d(x)^{-n}{\text
   d}V_F(x)$$
   and
$$ c(\alpha,r,R)= \int_{B(x_0,R)\setminus B(x_0,r)}d(x)^\alpha (\boldsymbol{\Delta}( \psi(x)d(x)^{-\gamma}))^2{\text d}V_F(x).$$ Clearly, $c(\alpha,r,R)$ is finite.   On the other hand,
\begin{eqnarray*}
  I_2(\varepsilon)&:=&\int_{M}{d(x)^{\alpha-4}}{\psi(x)^2u_\varepsilon(x)^2}{\text d}V_F(x)\\
   &\geq& \int_{B(x_0,r)\setminus  B(x_0,\varepsilon)}d(x)^{\alpha-4-2\gamma}{\text
d}V_F(x)\\&=&\tilde I(\varepsilon).
\end{eqnarray*}
By applying the layer cake representation and the volume comparison
(see Theorem \ref{comparison-laplace} (b)), we deduce that
\begin{eqnarray*}
\tilde I(\varepsilon)&=&\int_{B(x_0,r)\setminus
B(x_0,\varepsilon)}d(x)^{-n}{\text d}V_F(x)   = \int_{r^{-n}}^{\varepsilon^{-n}}{\rm Vol}_F(B(x_0,\rho^{-\frac{1}{n}})){\text d}\rho \\
   &\geq& \omega_n\int_{r^{-n}}^{\varepsilon^{-n}}\rho^{-1}{\text d}\rho\\
&=&n\omega_n(\ln r-\ln \varepsilon).
\end{eqnarray*}
In particular, $\lim_{\varepsilon\to 0^+}\tilde
I(\varepsilon)=+\infty.$ Therefore, it follows that
\begin{eqnarray*}
  \tilde C &\leq & \inf_{u\in
    C_{0,F,\alpha}^\infty(M)\setminus \{0\}}\frac{\int_{M}d(x)^\alpha (\boldsymbol{\Delta}u(x))^2{\text d}V_F(x)}{\int_{M}{d(x)^{\alpha-4}}{u(x)^2}{\text
    d}V_F(x)} \\
   &\leq& \lim_{\varepsilon\to
0^+}\frac{I_1(\varepsilon)}{I_2(\varepsilon)}\\
   &\leq&\lim_{\varepsilon\to
0^+}\frac{\gamma^2\left(\frac{n-\alpha}{2}+r\right)^2\tilde
   I(\varepsilon)+c(\alpha,r,R)}{\tilde
I(\varepsilon)}\\&=&\gamma^2\left(\frac{n-\alpha}{2}+r\right)^2.
\end{eqnarray*}
Since $r>0$ is arbitrary, we can take $r\to 0^+$,
which completes the proof of (\ref{best-const}). \hfill $\square$\\

Our second main result connects first to second order terms and it
can be stated as follows.

\begin{theorem}\label{main-2} {\rm (Rellich inequality II)} Let $(M,F)$ be an $n-$di\-men\-sional Finsler-Hadamard manifold with
${\bf S}=0$ and ${\bf K}\leq {c}\leq 0$, let $x_0\in M$ be fixed,
and choose any $\alpha\in \mathbb R$ such that $n-8+3\alpha>0$ and
$\alpha<2.$
 Then for every $u\in C_{0,F,\alpha}^\infty(M)$ we have
{ \begin{eqnarray*}
  \int_{M}d(x)^\alpha (\boldsymbol{\Delta}u(x))^2\text{\emph{d}}V_F(x) &\geq &
  \frac{(n-\alpha)^2}{4}\int_{M}{d(x)^{\alpha-2}}{F^*(x,Du(x))^2}\text{\emph{d}}V_F(x)\\&&+\frac{(n-4+\alpha)^2(n-\alpha)(n-1)}{8}\int_{M}{d(x)^{\alpha-4}}{\bf
D}_c(d(x)){u(x)^2}\text{\emph{d}}V_F(x).
\end{eqnarray*}}
Moreover, the constant $\frac{(n-\alpha)^2}{4}$ is sharp.
\end{theorem}

 {\it Proof.} We shall keep the notations and shall invoke some of the arguments from the proof of Theorem \ref{main-1}. Let $u\in C_{0,F,\alpha}^\infty(M)$. By applying the
 arithmetic-geometric mean inequality to the left hand side of
 (\ref{arith-mean}), it follows that
$$2\int_Md(x)^{\alpha-2}|u\boldsymbol{\Delta}(u(x))|{\text d}V_F(x)\leq \tilde C^{-\frac{1}{2}}\int_{M}d(x)^\alpha (\boldsymbol{\Delta}u(x))^2{\text d}V_F(x)
+\tilde C^{\frac{1}{2}}\int_{M}{d(x)^{\alpha-4}}{u(x)^2}{\text
d}V_F(x).$$ Combining this inequality with (\ref{mindjart-kell}), we
see that
\begin{eqnarray*}
  2\int_Md(x)^{\alpha-2}F^*(x,Du(x))^2{\text d}V_F(x) &\leq& \tilde
C^{-\frac{1}{2}}\int_{M}d(x)^\alpha
(\boldsymbol{\Delta}u(x))^2{\text d}V_F(x)\\&&+\left(\tilde C^{\frac{1}{2}}-2(2-\alpha)\gamma\right)\int_{M}{d(x)^{\alpha-4}}{u(x)^2}{\text d}V_F(x) \\
   &&-(2-\alpha)(n-1)\int_{M}{d(x)^{\alpha-4}}{\bf
D}_c(d(x)){u(x)^2}\text{\emph{d}}V_F(x).
\end{eqnarray*}
Since $\tilde
C^{\frac{1}{2}}-2(2-\alpha)\gamma=\frac{(n-8+3\alpha)\gamma}{2}>0$,
by applying Rellich inequality I to the second integrand on the
right hand side of the above inequality, a reorganization of the
expressions implies that
\begin{eqnarray*}
  2\int_Md(x)^{\alpha-2}F^*(x,Du(x))^2{\text d}V_F(x) &\leq& \frac{8}{(n-\alpha)^2}\int_{M}d(x)^\alpha
(\boldsymbol{\Delta}u(x))^2{\text d}V_F(x)\\
   &&-\frac{(n-4+\alpha)^2(n-1)}{n-\alpha}\int_{M}{d(x)^{\alpha-4}}{\bf
D}_c(d(x)){u(x)^2}\text{\emph{d}}V_F(x).
\end{eqnarray*}
Once we multiply this inequality by $\frac{(n-\alpha)^2}{8},$ we
obtain the Rellich inequality II.

It remains to prove that in Rellich inequality II the constant
$\frac{(n-\alpha)^2}{4}$ is sharp. By using the same functions as in
the proof of Theorem \ref{main-1}, it follows by
(\ref{tavolsag-derivalt}) that
\begin{eqnarray*}
  I_3(\varepsilon)&:=&\int_{M}{d(x)^{\alpha-2}}{F^*(x,D(\psi u_\varepsilon)(x))^2}{\text d}V_F(x)\\
   &\geq& \gamma^2\int_{B(x_0,r)\setminus  B(x_0,\varepsilon)}d(x)^{\alpha-4-2\gamma}{\text
d}V_F(x)\\&=&\gamma^2\tilde I(\varepsilon).
\end{eqnarray*}
The rest of the proof is similar as for Theorem \ref{main-1}.  \hfill $\square$\\

{\it Proof of Theorem \ref{Rellich-intro-1}.}  Take in Theorems
\ref{main-1} and \ref{main-2} the value $\alpha=0$. By considering
the
 continued fraction representation of the
function $\rho\mapsto \coth(\rho)$, one has
$$\rho\coth(\rho)-1\geq \frac{3\rho^2}{\pi^2+\rho^2},\ \forall
\rho>0,$$ and this concludes the proof.  \hfill $\square$

\section{Concluding remarks and questions}\label{sect-concluding}

\begin{remark}\rm [{\it Tour of Rellich inequalities}]
The technical hypothesis $n-8+3\alpha>0$ is indispensable in the
proof of Theorem \ref{main-2}. However, we believe an alternative
proof should eliminate this assumption. Interestingly, Rellich
inequalities I and II are {\it deducible from each other} via the
Hardy inequality once the assumption $n-8+3\alpha>0$ holds. First,
we have seen that the proof of Theorem \ref{main-2} is obtained from
the statement of Theorem \ref{main-1}. Conversely, by Rellich
inequality II and Hardy inequality (see relation
(\ref{hardy-majd-a-vegen-kell})), we obtain {
\begin{eqnarray*}
  \int_{M}d(x)^\alpha (\boldsymbol{\Delta}u(x))^2\text{\emph{d}}V_F(x) &\geq &
  \frac{(n-\alpha)^2}{4}\int_{M}{d(x)^{\alpha-2}}{F^*(x,Du(x))^2}\text{\emph{d}}V_F(x)\\&&+\frac{(n-4+\alpha)^2(n-\alpha)(n-1)}{8}\int_{M}{d(x)^{\alpha-4}}{\bf
D}_c(d(x)){u(x)^2}\text{\emph{d}}V_F(x)\\
&\geq&\frac{(n-\alpha)^2\gamma^2}{4}\int_{M}{d(x)^{\alpha-4}}{u(x)^2}\text{{d}}V_F(x)\\&&+\left[\frac{(n-4+\alpha)^2(n-\alpha)(n-1)}{8}+\frac{(n-\alpha)^2}{4}\gamma(n-1)\right]
 \times
\\&&\ \ \ \times \int_{M}{d(x)^{\alpha-4}}{\bf
D}_c(d(x)){u(x)^2}\text{{d}}V_F(x)\\&=&
\frac{(n-4+\alpha)^2(n-\alpha)^2}{16}\int_{M}{d(x)^{\alpha-4}}{u(x)^2}\text{\emph{d}}V_F(x)\\&&+\frac{(n-4+\alpha)(n-\alpha)(n-2)(n-1)}{4}\int_{M}{d(x)^{\alpha-4}}{\bf
D}_c(d(x)){u(x)^2}\text{\emph{d}}V_F(x),
\end{eqnarray*}}
which is precisely Rellich inequality I. In particular, the
Euclidean Rellich inequalities (\ref{standard-rellich-1}) and
(\ref{standard-rellich-2}) can be considered to be equivalent
whenever $n\geq 9.$
\end{remark}

\begin{remark}\rm \label{cut-off} 
[{\it Rigidity}]  For a generic Finsler manifold $(M,F)$  the
 vanishing of Green-deflection $G_F$ (where
the function $K_F$ appears)  played a crucial role in Rellich
inequalities. As we have already pointed out in Proposition
\ref{prop-eukl-green}, $K_F\equiv 0$ if and only if $(M,F)$ is
Riemannian.  On account of this characterization  we believe that
the {\it full} Rellich inequality holds, i.e., $$
\frac{(n-4+\alpha)^2(n-\alpha)^2}{16}=\inf_{u\in
    C_{0}^\infty(M)\setminus \{0\}}\frac{\int_{M}d(x)^\alpha (\boldsymbol{\Delta}u(x))^2{\text d}V_F(x)}{\int_{M}{d(x)^{\alpha-4}}{u(x)^2}{\text
    d}V_F(x)},$$
if and only if $(M,F)$ is Riemannian. Note that in Theorem
\ref{main-1} only the set of functions $C_{0,F,\alpha}^\infty(M)$ is
considered while the latter relation is formulated for the {\it
entire} space $C_{0}^\infty(M)$.
\end{remark}

\begin{remark}\rm [{\it Mean value property vs. $K_F\equiv 0$ on Minkowski spaces}]
Let $(M,F)=(\mathbb R^n,F)$ be a Minkowski space. Recently, Ferone
and Kawohl \cite[p. 252]{Fer-Kaw} proved the mean value property for
$\boldsymbol{\Delta}$-harmonics whenever
\begin{equation}\label{F-K-assump}
\frac{\langle a,b\rangle}{F(a)F^*(b)} =\langle \nabla F(a),\nabla
F^*(b)\rangle,\ \ \forall a,b\in \mathbb R^n\setminus\{0\}.
\end{equation}
Here, $\langle\cdot, \cdot\rangle$ denotes the usual inner product
on $\mathbb R^n.$ Interestingly, one can show that
(\ref{F-K-assump}) is {\it equivalent} to $K_F\equiv 0$, see
relation (\ref{j-csillag}). Therefore, according to Proposition
\ref{prop-eukl-green},  no proper non-Euclidean class of Minkowski
norms can be delimited in \cite{Fer-Kaw} to verify the mean value
property. In fact, we conjecture that the validity of the mean value
property of $\boldsymbol{\Delta}$-harmonics on a Minkowski space
$(\mathbb R^n,F)$ holds if and only if  $(\mathbb R^n,F)$ is
Euclidean. This problem will be studied in a forthcoming paper.
\end{remark}

\begin{remark}\rm [{\it Nonreversible Finsler manifolds}] In order
to avoid further technicalities, we focused our study only to
reversible Finsler manifolds. However, by employing suitable
modifications in the proofs, we can state Hardy and Rellich
inequalities on not necessarily reversible Finsler manifolds.
\end{remark}

\vspace{1cm}

\noindent {\bf Acknowledgment.} The research of A. Krist\'aly is
supported by J\'anos Bolyai Research Scholarship of the Hungarian
Academy of Sciences. Both authors were supported by the Slovenian
Research Agency grants P1-0292-0101 and J1-5435-0101.

\noindent {\normalsize {\sc  Alexandru Krist\'aly}\\ Institute of Applied Mathematics, \'Obuda University, 1034 Budapest, Hungary}\\
Email: {\textsf{alexandru.kristaly@yahoo.com}\\

\noindent {\normalsize \sc Du\v san Repov\v s}\\  {\normalsize
Faculty of Education, and Faculty of Mathematics and Physics,
University of
 Ljubljana, P.O.B. 2964, 1001
 Ljubljana,  Slovenia}\\
 Email: {\textsf{dusan.repovs@guest.arnes.si}}

\begin{thebibliography}{99}

\bibitem{AIM} P. L. Antonelli, R. S. Ingarden,  M. Matsumoto, The theory of
sprays and Finsler spaces with applications in physics and biology,
FTPH 58, Kluwer Academic Publishers, 1993.

 \bibitem{BCS} D.~Bao, S.~S.~Chern, Z.~Shen, {Introduction to
Riemann--Finsler Geometry,} Graduate Texts in Mathematics, 200,
Springer Verlag, Berlin, 2000.

\bibitem{BMBR} G. Bonanno, G. Molica Bisci, V. R\u adulescu, Multiple solutions of generalized Yamabe
equations on Riemannian manifolds and applications to Emden-Fowler
problems. Nonlinear Anal. Real World Appl. 12 (2011), no. 5,
2656--2665.

\bibitem{Carron-JMPA} G. Carron, In\'egalit\'es de Hardy sur les vari\'et\'es riemanniennes non-compactes.
J. Math. Pures Appl. (9) 76 (1997), no. 10, 883--891.

\bibitem{Fer-Kaw} V. Ferone,  B. Kawohl, Remarks on a Finsler-Laplace. Proc. Amer. Math. Soc.  137  (2009),  no. 1, 247--253.

\bibitem{GM} N. Ghoussoub,  A. Moradifam, Bessel pairs and optimal Hardy
and Hardy-Rellich inequalities. Math. Ann. 349 (2011), no. 1, 1--57.




\bibitem{KO-TAMS-2009} I. Kombe, M. \"Ozaydin, Improved Hardy and Rellich
inequalities on Riemannian manifolds. Trans. Amer. Math. Soc. 361
(2009), no. 12, 6191--6203.

\bibitem{KO-TAMS-2013} I. Kombe, M. \"Ozaydin, Hardy-Poincar\'e, Rellich and uncertainty principle inequalities on Riemannian manifolds. Trans. Amer. Math. Soc. 365 (2013), no. 10,
5035--5050.

\bibitem{Kristaly} A. Krist\'aly, Sharp uncertainty principles on Finsler manifolds: the effect of
curvature. Preprint, 2013. arXiv:1311.6418v2

\bibitem{Ohta-Sturm}  S. Ohta, K.-T. Sturm, Heat flow on Finsler manifolds. Comm. Pure Appl. Math. 62 (2009), no. 10,
1386--1433.

\bibitem{Shen-volume} Z. Shen, Volume comparison and its applications in Riemann-Finsler geometry. Adv. Math. 128 (1997), no. 2, 306--328.

\bibitem{Shen-monograph} Z. Shen, The non-linear Laplacian for Finsler manifolds. The theory of Finslerian Laplacians and applications, 187--198.
Mathematics and Its Applications, 459. Kluwer, Dordrecht, 1998.

\bibitem{TZ} A. Tertikas, N. B. Zographopoulos, Best constants in the Hardy-Rellich inequalities
and related improvements. Adv. Math. 209 (2007) 407--459.

\bibitem{Wu-Xin} B. Y. Wu, Y. L.  Xin, Comparison theorems in Finsler geometry and their applications. Math. Ann. 337 (2007), no. 1, 177--196.

\bibitem{YSK-Kinaiak} Q. Yang, D. Su, Y. Kong, Hardy inequalities on Riemannian manifolds with
negative curvature, Commun. Contemp. Math. (2013), in press. DOI:
10.1142/S0219199713500430

\bibitem{Zhao-Shen} W. Zhao, Y. Shen, A universal volume comparison theorem
for Finsler manifolds and related results. Canad. J. Math. 65
(2013), no. 6, 1401--1435.



\end{thebibliography}
\end{document}